\documentclass[12pt,twoside]{amsart}
\usepackage{amssymb}
\usepackage{amscd}

\title[Toric extremal contractions]
{Combinatorial descriptions\\
of toric extremal contractions} 
\author{Hiroshi Sato} 
\thanks{The author is partly supported by the 
Grant-in-Aid for JSPS Fellows, The Ministry of 
Education, Science, Sports and Culture, Japan.}
\address{Department of Mathematics\\ 
 Tokyo Institute of Technology, 2-12-1 Oh-Okayama, 
Meguro-ku, Tokyo 152-8551, Japan}
\email{hirosato@math.titech.ac.jp}

\newtheorem{thm}{Theorem}[section]
\newtheorem{lem}[thm]{Lemma}

\newtheorem{prop}[thm]{Proposition}
\newtheorem{lem-def}[thm]{Definition-Lemma}

\theoremstyle{definition}
\newtheorem{ex}[thm]{Example}
\newtheorem{defn}[thm]{Definition}
\newtheorem{rem}[thm]{Remark}
\newtheorem*{ack}{Acknowledgments}       
 
\newtheorem*{notation}{Notation}         

\newtheorem{say}[thm]{}
\begin{document}
\bibliographystyle{amsalpha+}

\begin{abstract}
In this paper, we give explicit combinatorial descriptions 
for toric extremal contractions under the relative setting, 
where varieties do not need to be complete. 
Fujino's completion theorem is the key to the main result. 
As applications, we can generalize some of Musta\c{t}\v{a}'s 
results related to Fujita's conjecture on toric varieties 
for the relative case. 
\end{abstract}

\maketitle
\tableofcontents

\section{Introduction}\label{intro}

The purpose of this paper is to give explicit combinatorial 
descriptions for extremal contractions from toric varieties. 
It is well-known that this problem was studied in \cite{reid} 
almost completely. However, the case where varieties are 
not complete seems not to be treated anywhere because of 
the difficulty of describing a non-complete fan. 
In this paper, we avoid this difficulty by using 
the notion of {\em extremal primitive relations}. 
Thanks to this 
notion, the descriptions are much simpler than Reid's. 
This paper is 
a kind of supplement to \cite{fujisato}. 

Fujino's completion theorem for toric morphisms 
in \cite{fujino2} is important. 
By using this theorem, promblems come down to the case where 
varieties are complete. 

The content of this paper is as follows: Section \ref{junbi} is 
a section for preparation. We review the local descriptions of 
toric extremal contractions when varieties are complete. 
We introduce the notion of extremal primitive relations. 
It is useful for describing toric extremal contractions. 
In Section \ref{relmori}, we give explicit descriptions for 
toric extremal contractions. 
Since we can apply Fujino's theorem, 
the problem becomes more elementary. 
In Section \ref{appl}, we give some generalizations of 
\cite{mustata} related to Fujita's conjecture 
as applications of the results in Section \ref{relmori}.

\begin{notation} 
Here, we summarize the some notation which we will 
use in this paper. 
\begin{itemize}
\item[(1)] Let $\Delta$ be a fan. 
For a cone $\sigma\in\Delta$, we denote 
the corresponding closed orbit by ${\rm V}(\sigma)$. 
\item[(2)] For a variety $X$ (resp. morphism $f:X\to Y$), 
we denote the Picard number (resp. relative Picard number) 
by $\rho(X)$ (resp. $\rho(X/Y)$).
\item[(3)] The symbol 
$\mathbb R_{\geq 0}$ denotes the set of non-negative real numbers. 
The other symbols are similar.
\end{itemize}
\end{notation}

\begin{ack} 
The author would like to thank Professor Osamu Fujino for 
introducing him to this problem 
and giving useful comments. 
He also wishes to thank Professor 
Natsuo Saito 
for advice and encouragement. 
\end{ack}

\section{Preliminaries}\label{junbi}

In this section, we review the toric Mori theory in \cite{reid}. 
For fundamental properties of the toric geometry, 
see \cite{fulton} and \cite{oda}. 
We recommend \cite{fujino}, \cite{fujisato} and \cite{ma} 
for understanding the toric Mori theory. 

Let $X=X_{\Delta}$ be a $\mathbb{Q}$-factorial 
toric $n$-fold 
associated to a fan $\Delta$ in $N={\mathbb Z}^n$ over an 
algebraically closed field. Let 
${\rm G}(\Delta)$ be the set of primitive generators 
of $1$-dimensional cones in $\Delta$, and put 
${\rm G}(\sigma):=\sigma\cap{\rm G}(\Delta)$ for a cone 
$\sigma\in\Delta$. Assume that $X$ is {\em complete}. 
The following notion is useful for describing 
extremal rays of $X$. See \cite{bat}, \cite{cas} and 
\cite{sato} more precisely.

\begin{defn}
A non-empty subset $P\subset{\rm G}(\Delta)$ is a 
{\em primitive collection} if $P$ does not generate 
any cone in $\Delta$, while every proper subset of $P$ 
generates a cone in $\Delta$.
\end{defn}

Let $f:X\to Y$ be a projective toric morphism. 
For an extremal ray $R\subset{\rm NE}(X/Y)$, let $w\in\Delta$ 
be the $(n-1)$-dimensional cone corresponding to $R$. 
By the completeness of $X$, there exist exactly two 
maximal cones $\sigma_1,\ \sigma_2\in\Delta$ 
such that $w\prec\sigma_1$ and $w\prec\sigma_2$. We put 
${\rm G}(w):=\{v_1,\ldots,v_{n-1}\}$, 
${\rm G}(\sigma_1)\setminus{\rm G}(w):=\{v_{n}\}$ and 
${\rm G}(\sigma_2)\setminus{\rm G}(w):=\{v_{n+1}\}$. 
After rearranging the elements in ${\rm G}(w)$, we obtain 
the equality
$$c_1v_1+\cdots+c_{\alpha}v_{\alpha}+
c_{\beta+1}v_{\beta+1}+\cdots+c_nv_n+c_{n+1}v_{n+1}=0,$$
where $c_1,\ldots,c_{\alpha}\in{\mathbb Z}_{<0}$, 
$c_{\beta+1},\ldots,c_{n+1}\in{\mathbb Z}_{>0}$, 
${\rm gcd}(c_1,\ldots,c_{n+1})=1$ and 
$\{v_{\beta+1},\ldots,v_{n+1}\}$ is a primitive collection. 
We call this relation 
an {\em extremal primitive relation}. Moreover, the following 
proposition holds.

\begin{prop}[{\cite[Lemma 3.1]{cas}} and 
{\cite[Theorem 4.10]{sato}}]\label{primlem1}
If $Q$ is a primitive collection of $\Delta$ such that
$Q\neq\{v_{\beta+1},\ldots,v_{n+1}\}$ and 
$Q\cap\{v_{\beta+1},\ldots,v_{n+1}\}\neq\emptyset$, 
then $(Q\setminus\{v_{\beta+1},\ldots,v_{n+1}\})\cup
\{v_1,\ldots,v_{\alpha}\}$ contains a primitive collection.
\end{prop}

\begin{rem}
In \cite{bat}, \cite{cas} and 
\cite{sato}, the arguments are developed under the assumption 
that varieties are {\em smooth}. However, 
combinatorial results of them like Proposition \ref{primlem1} 
hold for the simplicial case too. 
\end{rem}

Let 
$\sigma_i$ be the cone generated by 
$\{v_{\beta+1},\ldots,v_{n+1}\}\setminus\{v_i\}$ 
for $\beta+1\le i\le n+1$. 
By definition, we have $\sigma_i\in\Delta$. 
Put $w':={\mathbb R}_{\ge 0}v_1+\cdots
+{\mathbb R}_{\ge 0}v_{\alpha}\prec w$. 
Proposition \ref{primlem1} is equivalent to the following. 
We can confirm this equivalence as in the proof of 
{\cite[Theorem 4.10]{sato}}.

\begin{prop}\label{primlem2}
For any cone $\sigma\in\Delta$ such that $w'\prec\sigma$, 
put $\sigma=w'+\sigma'+\tau$, where $\sigma',\ \tau\in\Delta$, 
${\rm G}(\sigma')\subset\{v_{\beta+1},\ldots,v_{n+1}\}$ and 
${\rm G}(\tau)\cap\{v_1,\ldots,v_{\alpha},
v_{\beta+1},\ldots,v_{n+1}\}=\emptyset$. 
Then, $w'+\sigma_i+\tau\in\Delta$ for any 
$\beta+1\le i\le n+1$.
\end{prop}

By using Propositions \ref{primlem1} and \ref{primlem2}, 
we can 
construct the extremal contraction $\varphi_R:X\to W$ 
associated to $R$ (see Section \ref{relmori}). 

\begin{say}[Extremal primitive relations]
Let $f:X\to Y$, $R$ and $w$ as above. We rewrite the extremal 
primitive relations with respect to $R$ as
$$a_1x_1+\cdots+a_lx_l=b_1y_1+\cdots+b_my_m,$$
where $\{x_1,\ldots,x_l,y_1,\ldots,y_m\}\subset{\rm G}(\Delta)$ 
and $a_1,\ldots,a_l,b_1,\ldots,b_m\in{\mathbb Z}_{>0}$. 
Namely, $l=n-\beta+1,$ $m=\alpha,$ $a_1x_1+\cdots+a_lx_l=
c_{\beta+1}v_{\beta+1}+\cdots+c_{n+1}v_{n+1}$ and 
$b_1y_1+\cdots+b_my_m=-(c_1v_1+\cdots+c_{\alpha}v_{\alpha})$. 
The following is obvious.

\begin{prop}\label{hotei}
For a torus invariant prime divisor $D$ which corresponds to 
$v\in{\rm G}(\Delta)$ and a curve $C$ which spans $R$, 
the following holds$:$
$$
(D\cdot C)
\begin{cases}
>0 \\
=0 \\
<0 \\
\end{cases}
\quad {\text{if and only if}} \quad
\begin{cases}
v\in\{x_1,\ldots,x_l\},\\
v\not\in\{x_1,\ldots,x_l,y_1,\ldots,y_m\},\\
v\in\{y_1,\ldots,y_m\}.\\
\end{cases}
$$ 
\end{prop}

We can completely recover $\Delta$ from three data: 
$N$, the primitive collections and the extremal 
primitive relations. This follows from the well-known 
isomorphism
$$
{\rm A}_1(X)\otimes{\mathbb Q}\cong
\left\{(c_x)_{x\in{\rm G}(\Delta)}
\in{\mathbb Q}^{{\rm G}(\Delta)}\,\left|\,
\sum_{x\in{\rm G}(\Delta)}c_xx=0\right.\right\},
$$
where ${\rm A}_1(X)$ is the Chow group of $1$-cycles. 
In particular, if ${\rm G}(\Delta)$ generates $N$, 
then we can recover $\Delta$ from 
the primitive collections and the extremal 
primitive relations. If $\rho(X)=1$, 
then this condition is equivalent to the one that 
$X$ is a weighted projective space. This is well-known.

\begin{ex}
Let $X=X_{\Delta}$ be a ${\mathbb Q}$-factorial terminal 
toric ${\mathbb Q}$-Fano $3$-fold with 
Picard number $1$ whose extremal primitive relation is 
$x_1+x_2+x_3+x_4=0$. Then, there exist exactly 
two possibilities for such 
$\Delta$ (see \cite{kaspr}). 
\end{ex}

\end{say}

\section{Relative toric Mori theory}\label{relmori}

In this section, we deal with the relative toric Mori theory 
from the combinatorial viewpoint. 
We remark that varieties are {\em not necessarily complete}. 
For the general theory, see 
\cite{fujisato}. 

Let $f:X=X_{\Delta}\to Y$ be 
a projective surjective toric morphism 
with $\dim X=n$. 
We assume that $X$ is ${\mathbb Q}$-factorial. 
For an extremal ray $R\subset{\rm NE}(X/Y)$, 
let $\varphi:=\varphi_R:X\to W$ be
the associated extremal contraction. The following Fujino's 
theorem in \cite{fujino2} is the key to the main result 
of this section. 

\begin{thm}[{\cite[Theorem 2.10]{fujino2}}]\label{godfujino}
There exist equivariant completions 
$$
\begin{matrix} 
\overline{X} & \stackrel{\overline{\varphi}}{\longrightarrow} 
& \ \overline{W} \\
\ \ \ \ \searrow & \ &  {\swarrow}\ \ \ \ \\
 \  & \overline{Y} &  
\end{matrix}
$$
of toric morphisms
$$
\begin{matrix} 
X & \stackrel{\varphi}{\longrightarrow} 
& \ W, \\
\ \ \  \ \searrow & \ &  {\swarrow}\ \ \ \ \\
 \  & Y &  
\end{matrix}
$$
where 
\begin{enumerate}
\item[(1)] $\overline{X}$, $\overline{Y}$ and $\overline{W}$ are 
equivariant completions of $X$, $Y$ and $W$, respectively, 
\item[(2)] $\overline{X}$ is ${\mathbb Q}$-factorial, 
\item[(3)] $\overline{\varphi}:\overline{X}\to\overline{W}$, 
$\overline{X}\to\overline{Y}$ and $\overline{W}\to\overline{Y}$ 
are projective and 
\item[(4)] $\rho(\overline{X}/\overline{W})=1$.
\end{enumerate}

\end{thm}

Thus, we fix an equivariant completion 
$\overline{\varphi}:\overline{X}\to\overline{W}$ of $\varphi$ 
as in Theorem \ref{godfujino}. 
Let $\overline{\Delta}$ be the fan associated to $\overline{X}$ 
and $w\in\overline{\Delta}$ an $(n-1)$-dimensional cone 
such that $\overline{\varphi}({\rm V}(w))$ is a point. Then, 
as in Section \ref{junbi}, we have an extremal primitive 
relation 
$a_1x_1+\cdots+a_lx_l=b_1y_1+\cdots+b_my_m$ 
for $w$, where $x_1,\ldots,x_l,
y_1,\ldots,y_n\in{\rm G}(\overline{\Delta})$. 
We define $w'\in\overline{\Delta}$ 
and $\sigma_i\in\overline{\Delta}$ for $1\le i\le l$ 
similarly as in Section \ref{junbi}.

\begin{lem}
$w'+\sigma_i\in\Delta$ for any $1\le i\le l$. 
\end{lem}
\begin{proof}
Obviously, $x_1,\ldots,x_l,
y_1,\ldots,y_m$ are contained in the inverse image of 
a cone in $\Delta_Y$. 
So, by the properness of $f$, 
we complete the proof. 
\end{proof}

Thus, we call 
$a_1x_1+\cdots+a_lx_l=b_1y_1+\cdots+b_my_m$ 
an {\em extremal primitive relation} in this case, too.

Proposition \ref{primlem2} also holds for this case. 
This immediately follows from the properness of $\varphi$ 
and Proposition \ref{primlem2}. Though the setting is 
distinct from the one in Proposition \ref{primlem2}, 
the statement is completely similar. We repeat it. 

\begin{thm}\label{primlem3}
For any cone $\sigma\in\Delta$ such that $w'\prec\sigma$, 
put $\sigma=w'+\sigma'+\tau$, where $\sigma'$ and $\tau$ 
are as in Proposition $\ref{primlem2}$. 
Then, $w'+\sigma_i+\tau\in\Delta$ for any 
$1\le i\le l$. 
\end{thm}

In the remaining part of this section, we give an explicit 
combinatorial description for $\varphi:X\to W$ by using 
Theorem \ref{primlem3}. 

\begin{say}[Fano contractions]
Suppose that $\varphi$ is a {\em Fano} contraction, 
that is, $\dim X>\dim W$. 
This condition is equivalent to $m=0$. 
Put $N':=N/({\mathbb Z}x_1+\cdots+{\mathbb Z}x_{l})$. 
Then, we obtain the fan $\Delta_W$ in $N'$ associated to $W$ 
by sending the cones in $\Delta$ through $N\to N'$. 
The general fiber $F$ of $\varphi$ is the complete 
toric variety of Picard number $1$ 
whose extremal primitive relation is 
$a_1x_1+\cdots+a_lx_l=0$. In particular, 
if $\{x_1,\ldots,x_l\}$ generates the lattice, then 
$F\cong{\mathbb P}(a_1,\ldots,a_l)$ 
(see Section \ref{junbi}). 
\end{say}

\begin{say}[Birational contractions]
Suppose that $\varphi$ is {\em birational}. 
This condition is equivalent to $m>0$. 
Put 
$$\widetilde{w}:={\mathbb R}_{\ge 0}y_1+\cdots+
{\mathbb R}_{\ge 0}y_{m}+{\mathbb R}_{\ge 0}x_1
+\cdots+{\mathbb R}_{\ge 0}x_l.$$
We remark that if $m=1$, then 
$y_1\in{\mathbb R}_{\ge 0}x_{1}
+\cdots+{\mathbb R}_{\ge 0}x_{l}$. 
For $w'\prec\sigma\in\Delta$, put 
$\sigma'$ and $\tau$ as in Theorem $\ref{primlem3}$. 
The fan $\Delta_W$ in $N$ associated to $W$ is as follows:
$$\Delta_W=\left(\Delta\setminus
\left\{\sigma\in\Delta\,\left|\,w'\prec\sigma
\right.\right\}\right)\cup\left\{\widetilde{w}+
\tau\,\left|\,w'\prec\sigma=w'+\sigma'+\tau\in
\Delta\right.\right\}.$$
The exceptional locus $A$ of $\varphi$ is ${\rm V}(w')$, 
while $B:=\varphi(A)={\rm V}(\widetilde{w})$. We have 
${\rm codim}\,A=m$ and $\dim B=n-l-m+1$. 
We note that $A$ is irreducible. This does not 
necessarily hold for non-${\mathbb Q}$-factorial 
varieties (see \cite[Example 4.1]{fujino2}). 

\end{say}

\begin{say}[Flips, flops and anti-flips]
Suppose that $\varphi$ is a {\em small} contraction, 
that is, ${\rm codim}\,A\ge 2$. 
This condition is equivalent to $m\ge 2$. 
We use the same notation as above. Put 
$$w^+:={\mathbb R}_{\ge 0}x_1+\cdots+
{\mathbb R}_{\ge 0}x_l.$$
Then, we obtain a fan $\Delta^+$ in $N$ by star-subdividing 
$\Delta_W$ along $w^+$. Let $X^+$ be the associated 
toric variety. Then, the diagram
$$
\begin{matrix} 
X & \dashrightarrow & \ X^+ \\
{\ \ \ \ \searrow} & \ &  {\swarrow}\ \ \ \ \\
 \ & W &  
\end{matrix}
$$
is 
$$
\begin{cases}
\text{a flip for } \varphi \\
\text{a flop for } \varphi \\
\text{an anti-flip for } \varphi \\
\end{cases}
\quad {\text{if}} \quad 
a_1+\cdots+a_l-(b_1+\cdots+b_m)
\begin{cases}
>0,\\
=0,\\
<0.\\
\end{cases}
$$
We remark that 
$b_1y_1+\cdots+b_my_m=a_1x_1+\cdots+a_lx_l$ 
is an extremal primitive relation of $\Delta^+$. 

\end{say}

\section{Applications}\label{appl}

As applications of the results in the previous section, 
we can generalize some results in \cite{mustata}
 (cf. \cite[Remark 3.3]{fujino} and 
\cite[Theorem 3.9]{fujino2}). 
We use the same notation as in Sections \ref{junbi} and 
\ref{relmori}. In addition, for the case of a Fano 
contraction, put $A:=X$ and $B:=W$.

Let $f:X\longrightarrow Y$ be a projective surjective 
toric morphism as in the previous sections. We assume that 
$X$ is {\em smooth}. For an extremal ray of ${\rm NE}(X/Y)$, 
let $\varphi_R:X\longrightarrow W$ be the extremal 
contraction with respect to $R$. 
Then, the general fiber of 
$\varphi_R:A\longrightarrow B$ is a projective space, 
since $X$ is smooth. So, let $C_R$ be 
a line in a general fiber 
of $A\longrightarrow B$. 
The following is obvious.

\begin{lem}\label{mus02}
$C_R$ spans $R$, and 
for any torus invariant prime divisor $D$ on $X$, 
we have $(D\cdot C_R)\le 1$.
\end{lem}

\begin{rem}
If $R$ contains a numerical class of torus 
invariant curves, then the local description of 
the extremal contraction $\varphi_R$ 
coincides with Reid's. Therefore, 
if $R$ contains a numerical class of 
torus invariant curves, then we can make $C_R$ 
torus invariant. This is obvious by Reid's description of 
$\varphi_R$. 
\end{rem}

\begin{prop}\label{mus03}
Let $f:X\longrightarrow Y$ be as above, $L$ an 
line bundle on $X$ and $l$ a positive integer. 
Assume that $X$ is smooth. If $(L\cdot C_R)\ge l$ for 
every 
extremal ray $R$ of ${\rm NE}(X/Y)$, 
then for every torus invariant prime divisor $D$ on $X$, 
we have $(L(-D)\cdot C_R)\geq l-1$. 
In particular, if $L$ is $f$-ample, then $L(-D)$ is $f$-free. 
\end{prop}

\begin{proof}
The first part is obvious by Lemma \ref{mus02}. The last part 
follows from the equivalence of $f$-freeness and $f$-nefness 
on toric varieties. 
\end{proof}

Proposition \ref{mus03} is a relative version of 
\cite[Proposition 4.3 and Lemma 4.4]{mustata}. 

\begin{prop}\label{mus04}
Let $f:X\longrightarrow Y$ be as above and $L$ an 
$f$-ample line bundle on $X$. Assume that $X$ is smooth. 
For two distinct 
torus invariant prime divisors $D_1$ and $D_2$ on $X$, 
let $v_1$ and $v_2$ be the corresponding elements in 
${\rm G}(\Delta)$, respectively. 
Then, $L(-D_1-D_2)$ is not $f$-free if and only if 
there exists an extremal ray $R$ of ${\rm NE}(X/Y)$ whose 
extremal primitive relation is 
$a_1x_1+\cdots+a_lx_l=b_1y_1+\cdots+b_my_m$ such that 
$\{v_1,v_2\}\subset\{x_1,\ldots,x_l\}$ and $(L\cdot C_R)=1$.
\end{prop}

\begin{proof}
If there exists an extremal ray $R$ as above, then 
$(L(-D_1-D_2)\cdot C_R)=-1$. Therefore, $L(-D_1-D_2)$ is 
not $f$-nef, that is, not $f$-free. 

So, suppose that $L(-D_1-D_2)$ is not $f$-nef. There exists 
an extremal ray $R$ of ${\rm NE}(X/Y)$ such that 
$(L(-D_1-D_2)\cdot C_R)<0$. Since $(L\cdot C_R)\ge 1$, 
$(D_1\cdot C_R)\le 1$ and $(D_2\cdot C_R)\le 1$, we have 
$(L\cdot C_R)=(D_1\cdot C_R)=(D_2\cdot C_R)=1$. This completes 
the proof by Proposition \ref{hotei}. 
\end{proof}

Proposition \ref{mus04} is the generalization of 
\cite[Proposition 4.5]{mustata}, while the following is the one 
of \cite[Proposition 4.6]{mustata}. The proof is similar. 

\begin{prop}\label{mus05}
Let $f:X\longrightarrow Y$ be as above and $L$ an 
$f$-ample line bundle on $X$. Assume that $X$ is smooth. 
For a torus invariant prime divisor $D$ on $X$, 
let $v$ be the corresponding element in 
${\rm G}(\Delta)$. 
Then, $L(-D)$ is not $f$-ample if and only if 
there exists an extremal ray $R$ of ${\rm NE}(X/Y)$ whose 
extremal primitive relation is 
$a_1x_1+\cdots+a_lx_l=b_1y_1+\cdots+b_my_m$ such that 
$v\in\{x_1,\ldots,x_l\}$ and $(L\cdot C_R)=1$.
\end{prop}

\end{document}